\documentclass{article}

\usepackage{amssymb, latexsym}

\usepackage{graphicx}

\usepackage{amsmath}
\usepackage{color}

\newtheorem{theorem}{Theorem}[section]

\newtheorem{proposition}[theorem]{Proposition}

\newtheorem{remark}[theorem]{Remark}

\begin{document}
\title{\large\bf On the Shannon entropy power on Riemannian manifolds and Ricci flow}
\author{\ \ Songzi Li, Xiang-Dong Li
\thanks{Research supported by NSFC No. 11771430, Key Laboratory RCSDS, CAS, No. 2008DP173182, and a
Hundred Talents Project of AMSS, CAS.} }

\maketitle

\begin{minipage}{120mm}
{\bf Abstract}
In this paper, we prove the  concavity  of the Shannon entropy power  for the heat equation associated with the Laplacian or  the Witten Laplacian on complete Riemannian manifolds with suitable curvature-dimension condition and on compact super Ricci flows. Under suitable  curvature-dimension condition, we prove that the rigidity models of the Shannon entropy power are Einstein or quasi Einstein manifolds with Hessian solitons. Moreover, we prove  the convexity of the Shannon entropy power for the conjugate heat equation introduced by G. Perelman on Ricci flow and that  the corresponding rigidity models are the shrinking Ricci solitons. As an application, we prove the entropy isoperimetric inequality on complete Riemannian manifolds with non-negative (Bakry-Emery) Ricci curvature  and the maximal volume growth condition. 
\end{minipage}




\section{Introduction}

In his 1948 seminal paper \cite{Sh48},  Shannon introduced the notion of entropy power for  continuous random vectors and discovered the {\it Entropy Power Inequality}  (EPI). More precisely, let $X$ be an $n$-dimensional  continuous random vector with probability distribution $f(x)dx$, and 
\begin{eqnarray*}
H(X)=H(f)=-\int_{\mathbb{R}^n} f(x)\log f(x)dx \label{HX}
\end{eqnarray*}
 be the Boltzmann-Shannon differential entropy of $X$ or $f$.  The Shannon entropy power of $X$ or $f$ is defined as follows
\begin{eqnarray*}
N(X)=N(f)=e^{{2\over n}H(X)}.\label{NX}
\end{eqnarray*}
The Entropy Power Inequality (EPI) can be stated as follows: Let $X$ and $Y$ be two independent continuous random vectors with values in $\mathbb{R}^n$.  Then
\begin{eqnarray}
N(X+Y)\geq N(X)+N(Y). \label{NXY1}
\end{eqnarray}
Equivalently, for any two probability density functions $f$ and $g$ on $\mathbb{R}^n$, we have 
\begin{eqnarray}
N(f*g)\geq N(f)+N(g), \label{NXY3}
\end{eqnarray}
where $f*g$ denotes the convolution of $f$ and $g$. For the first complete proof of the Entropy Power inequality $(\ref{NXY1})$ or $(\ref{NXY3})$, see 
Stam  \cite{Stam} and Blachman \cite{Bla}. See also \cite{CT}.  

%
%
%
%
%

In \cite{Cost}, Costa proved the concavity of the Shannon entropy power along the heat equation on $\mathbb{R}^n$. More precisely, let  $u(x, t)$ be  the unique solution to the heat equation on $\mathbb{R}^n$
\begin{eqnarray*}
\partial_t u=\Delta u, \ \ \ u(x, 0)=f(x). \label{heat1}
\end{eqnarray*}
Let 
\begin{eqnarray*}
H(u(t))=-\int_{\mathbb{R}} u\log u dx \label{Hu}
\end{eqnarray*}
be the differential entropy associated to the heat distribution $u(x, t)dx$ at time $t$, and 
\begin{eqnarray*}
N(u(t))=e^{{2\over n}H(u(t))}. \label{Nu}
\end{eqnarray*}
be the Shannon entropy power of  $u(x, t)dx$ at time $t$. 
Then
\begin{eqnarray}
{d^2\over dt^2}N(u(t))\leq 0. \label{N2}
\end{eqnarray}
%
%

Using an argument based on the Blachman-Stam inequality \cite{Bla}, the original proof of the Entropy Power Concave Inequality (EPCI) $(\ref{N2})$ has been simplified by Dembo et al. \cite{Dem1, Dem2} and Villani \cite{V00}. 
In \cite{V00},  Villani pointed out the possibility of extending the Entropy Power Concave Inequality (EPCI) $(\ref{N2})$ to Riemannian manifolds with non-negative Ricci curvature using the $\Gamma_2$-calculation. In \cite{ST}, the concavity property of the entropy power has been also extended to Renyi entropy power, when evaluated along the solution to a nonlinear diffusion equation,  the so-called porous medium equation.

The purpose of this paper is to prove the  concavity  of the Shannon entropy power  for the heat equation associated with the Laplacian or  the Witten Laplacian on complete Riemannian manifolds with suitable curvature-dimension condition and on compact super Ricci flows. Under suitable curvature-dimension condition, we prove that the rigidity models of the Shannon entropy power are Einstein or quasi Einstein manifolds with Hessian solitons. Moreover, we prove the convexity of the Shannon entropy power for the conjugate heat equation introduced by G. Perelman on Ricci flow and the corresponding rigidity models are the shrinking Ricci solitons. As an application, we prove the entropy isoperimetric inequality on complete Riemannian manifolds with non-negative (Bakry-Emery) Ricci curvature  and the maximal volume growth condition. 
This might lead us to recognize  the importance of the {\it information-theoretic  approach} in the future study of geometric analysis  and Ricci flow. 

\section{Notation and main results}

Let $(M, g)$
be an $n$-dimensional complete Riemannian manifold, $\phi\in C^2(M)$ and $d\mu=e^{-\phi}dv$, where $v$ is the Riemannian
volume measure on $(M, g)$.
The Witten Laplacian acting on smooth functions is defined by
$$L = \Delta - \nabla\phi\cdot\nabla.$$
For any $u, v \in C_0^\infty(M)$, the integration by parts formula holds
\begin{eqnarray*}
\int_M \langle \nabla u, \nabla v\rangle d\mu=-\int_M L u vd\mu= - \int_M u L v d\mu.
\end{eqnarray*}
Thus, $L$ is the infinitesimal generator of the Dirichlet form
\begin{eqnarray*}
\mathcal{E}(u, v)=\int_M \langle\nabla u, \nabla v\rangle d\mu, \ \ \ \ u, v\in C_0^\infty(M).
\end{eqnarray*}
By It\^o's theory, the Stratonovich SDE on $M$
\begin{eqnarray*}
dX_t =\sqrt{2}U_t\circ dW_t-\nabla\phi(X_t)dt,\ \ \ \ \ \nabla_{\circ dX_t} U_t=0,
\end{eqnarray*}
where $U_t$ is the stochastic parallel transport along the trajectory of $X_t$, with initial data $X_0=x$ and $U_0={\rm Id}_{T_xM}$, defines a diffusion process $X_t$ on $M$ with infinitesimal generator $L$.
Moreover,  the  transition probability density of 
the $L$-diffusion process $X_t$ with respect to $\mu$, i.e., the heat kernel  $p_t(x, y)$   of the Witten Laplacian $L$, is the fundamental solution to the heat equation
\begin{eqnarray}
\partial_t u=Lu.  \label{heq}
\end{eqnarray}

In \cite{BE}, Bakry and Emery proved the generalized Bochner formula 
\begin{eqnarray}
L|\nabla u|^2-2\langle \nabla u, \nabla L u\rangle=2\|\nabla^2
u\|_{\rm HS}^2+2Ric(L)(\nabla u, \nabla u), \label{BWF}
\end{eqnarray}
where $u\in C^2(M)$, $\nabla^2 u$ denotes the Hessian of $u$, $\|\nabla^2 u\|_{\rm HS}$ is its Hilbert-Schmidt norm, and 
$$Ric(L)= Ric + \nabla^2\phi$$
is now called the (infinite dimensional) Bakry-Emery Ricci curvature associated with the Witten Laplacian $L$.  For $m\in [n, \infty)$,  the $m$-dimensional Bakry-Emery Ricci curvature associated with the Witten Laplacian $L$  is defined by
$$
Ric_{m, n}(L) = Ric + \nabla^2\phi - {\nabla\phi\otimes\nabla\phi\over m-n}.
$$
In view of this, we have
\begin{eqnarray*}
L|\nabla u|^2-2\langle \nabla u, \nabla L u\rangle \geq {2|Lu|^2\over m}+2Ric_{m, n}(L)(\nabla u, \nabla u).
\end{eqnarray*}
Here we make a convention that $m=n$ if and only if $\phi$ is a constant. By definition, we have
$$Ric(L)=Ric_{\infty, n}(L).$$
Following \cite{BE}, we say that $(M, g, \phi)$ satisfies the curvature-dimension 
$CD(K, m)$-condition for a constant $K\in \mathbb{R}$ and $m\in [n, \infty]$ if and only if
$$Ric_{m, n}(L)\geq Kg.$$
Note that, when $m=n$, $\phi=0$, we have $L=\Delta$ is the usual Laplacian on $(M, g)$, and the $CD(K, n)$-condition holds if and only if the Ricci curvature on $(M, g)$ is bounded from below by $K$, i.e., 
$$Ric\geq Kg.$$

In the case of  Riemannian manifolds with a family of time dependent metrics and potentials, we call $(M, g(t), \phi(t), t\in [0, T])$ a $(K, m)$-super Ricci flow if the metric $g(t)$ and the potential function $\phi(t)$ satisfy 
\begin{eqnarray}
{1\over 2}{\partial g\over \partial t}+Ric_{m, n}(L)\geq Kg, \label{KmsRF}
\end{eqnarray}
where $$L=\Delta_{g(t)}-\nabla_{g(t)}\phi(t)\cdot\nabla_{g(t)}$$ is the time dependent Witten Laplacian on $(M, g(t), \phi(t), t\in [0, T])$, and $K\in \mathbb{R}$ is a constant. When $m=\infty$, i.e., 
if the metric $g(t)$ and the potential function $\phi(t)$ satisfy the following inequality
\begin{eqnarray*}
{1\over 2}{\partial g\over \partial t}+Ric(L)\geq Kg,
\end{eqnarray*}
we call $(M, g(t), \phi(t), t\in [0, T])$   a $(K, \infty)$-super Ricci flow or a $K$-super Perelman Ricci flow. Indeed the $(K, \infty)$-Ricci flow (called also the $K$-Perelman Ricci flow)  
\begin{eqnarray*}
{1\over 2}{\partial g\over \partial t}+Ric(L)=Kg
\end{eqnarray*}
is a natural extension of the modified Ricci flow $\partial_t g=-2Ric(L)$ introduced by Perelman \cite{P1}  as the gradient flow of $\mathcal{F}(g, \phi)=\int_M (R+|\nabla \phi|^2)e^{-\phi}dv$ on $\mathcal{M}\times C^\infty(M)$ 
under the constraint condition that the measure $d\mu=e^{-\phi}dv$ is preserved. For the study of the Li-Yau or Hamilton differential Harnack inequalities, $W$-entropy formulas and related functional inequalities on $(K, m)$ or $(K, \infty)$-super Ricci flows, see \cite{LL15, LL-AJM, LL17,  LL19, LL18a, LL18b} and references therein. For super Ricci flows on metric and measure spaces, see \cite{Sturm18} and references therein. 

\medskip

Now we state the main results of this paper.

\begin{theorem} \label{thm1} Let $M$ be a complete Riemannian manifold with $Ric\geq Kg$ for some $K\in \mathbb{R}$.  Let  $u$ be a positive solution to the heat equation
$\partial_t u=\Delta u$. Let
\begin{eqnarray*}
H(u(t))=-\int_M u\log ud\mu, \ \ \ \ N(u(t))=e^{{2\over n}H(u(t))}.
\end{eqnarray*}
Then the Entropy Differential Inequality ${\rm EDI}(K, n)$ holds  on $(0, \infty)$ 
\begin{eqnarray}
H''+{2H'^2\over n}+2KH'\leq 0   \label{HNK1}
\end{eqnarray}
with the initial boundary condition $\lim\limits_{t\rightarrow 0+}tH'(u(t))\leq {n\over 2}$.
Equivalently,  the Entropy Power Concavity Inequality ${\rm EPCI}(K, n)$ holds on $(0, \infty)$ 
\begin{eqnarray}
{d^2 N\over dt^2}\leq -2K {dN\over dt}.  \label{NKN1}
\end{eqnarray}
Moreover, the equality in $(\ref{HNK1})$ or $(\ref{NKN1})$ holds on $(0, T]$ for some $T>0$  if and only if $(M, g)$ is a Einstein manifold with Hessian soliton
\begin{eqnarray}
Ric=Kg, \ \ \ \ \nabla^2 f={H'_{K}(t)\over n}g,
\end{eqnarray}
where $H'_K$ is the solution to the entropy differential  equation 
\begin{eqnarray*}
H_K''+{2\over n}H_K'^2+2KH'_K=0, \ \ \ \lim\limits_{t\rightarrow 0+} tH_K'(t)\leq {n\over 2}.
\end{eqnarray*}
In particular, if $Ric\geq 0$,  then $N(u(t))$ is concave on $(0, \infty)$, i.e., 
\begin{eqnarray}
{d^2N\over dt^2} \leq 0.  \label{NCI1}
\end{eqnarray}
Moreover, the equality in $(\ref{NCI1})$ holds on $(0, T]$ for some $T>0$  if and only if $M$ is isometric to $\mathbb{R}^n$, and 
\begin{eqnarray*}
f(x, t)={\|x\|^2\over 4t},  \ \ \ \ \forall x\in \mathbb{R}^n, t>0.
\end{eqnarray*}
That is to say, the equality in $(\ref{NCI1})$ holds if and only if $M$ is the Euclidean space $\mathbb{R}^n$ with the Gaussian Ricci soliton. 
\end{theorem}

\begin{theorem}\label{thm2} Let $(M, g)$ be a complete Riemannian manifold with $Ric_{m, n}(L)\geq Kg$ for some constants $m\geq n$ and $K\in \mathbb{R}$. Let  $u$ be a positive solution to the heat equation
$\partial_t u=Lu$. Let
\begin{eqnarray*}
H(u(t))=-\int_M u\log ud\mu, \ \ \ \ N(u(t))=e^{{2\over m}H(u(t))}.
\end{eqnarray*}Then
the Entropy Differential Inequality EDI$(K, m)$ holds
\begin{eqnarray}
H''+{2H'^2\over m}+2KH'\leq 0,\label{EDI}
\end{eqnarray}
with the initial boundary condition $\lim\limits_{t\rightarrow 0+}tH'(u(t))\leq {m\over 2}$. 
Equivalently, the Entropy Power  Concavity Inequality EPCI$(K, m)$ holds  on $(0, \infty)$
\begin{eqnarray}
{d^2 N \over dt^2}\leq-2K {dN \over dt}.  \label{ECI1}
\end{eqnarray}
Moreover, , the equality in $(\ref{EDI})$ or $(\ref{NCI1})$ holds on $(0, T]$ for some $T>0$  if and only if $(M, g)$ is a quasi-Einstein manifold with Hessian soliton
\begin{eqnarray*}
Ric_{m, n}(L)=Kg, \ \ \ \ \nabla^2 f={H'_{K}(t)\over m}g,
\end{eqnarray*}
where $H'_K$ is the solution to the entropy differential  equation 
\begin{eqnarray*}
H_K''+{2\over m}H_K'^2+2KH'_K=0, \ \ \ \lim\limits_{t\rightarrow 0+} tH_K'(t)\leq {m\over 2}.
\end{eqnarray*}
In particular, if $Ric_{m, n}(L)\geq 0$,  then $N(u(t))$ is concave  on $(0, \infty)$, i.e., 
\begin{eqnarray}
{d^2N\over dt^2} \leq 0.  \label{NCI2}
\end{eqnarray}
Moreover, the equality in $(\ref{NCI2})$ holds on $(0, T]$ for some $T>0$ if and only if $M$ is isometric to $\mathbb{R}^n$,  $\phi$ is a constant, $m=n$, and 
\begin{eqnarray*}
f(x, t)={\|x\|^2\over 4t} \ \ \ \ \forall x\in \mathbb{R}^n, t>0.
\end{eqnarray*}
That is to say, the equality in $(\ref{NCI2})$ holds on $(0, T]$ for some $T>0$ if and only if $M$ is Euclidean space with the Gaussian Ricci soliton. 
\end{theorem}

The following result extends Theorem \ref{thm1}  to compact $(K, m)$-super Ricci flows.

\begin{theorem}\label{thm3}  Let  $(M, g(t), \phi(t), t\in [0, T])$ be a compact $(K, m)$-super  Ricci flow in the sense that
\begin{eqnarray*}
 {1\over 2}{\partial g\over \partial t}+Ric_{m, n}(L)\geq Kg,\ \ \  {\partial \phi\over \partial t}={1\over 2}{\rm Tr}\left({\partial g\over \partial t}\right).
\end{eqnarray*}
then the Entropy Differential inequality $(\ref{EDI})$ and the Entropy  Power Concavity Inequality  $(\ref{ECI1})$ hold on $(0, T]$. 
Moreover, the equality in $(\ref{EDI})$ or $(\ref{ECI1})$  holds if and only if $(M, g(t), \phi(t))$ is a $(K, m)$-Ricci flow
\begin{eqnarray*}
 {1\over 2}{\partial g\over \partial t}+Ric_{m, n}(L)=Kg,\ \ \ {\partial \phi\over \partial t}={1\over 2}{\rm Tr}\left({\partial g\over \partial t}\right)=0,
\end{eqnarray*}
and $f$ is a Hessian soliton 
\begin{eqnarray*}
\nabla^2 f={H'_K(t)\over m}g,
\end{eqnarray*}
where $H'_K$ is the solution to the entropy differential  equation 
\begin{eqnarray*}
H_K''+{2\over m}H_K'^2+2KH'_K=0, \ \ \ \lim\limits_{t\rightarrow 0+} tH_K'(t)\leq {m\over 2}.
\end{eqnarray*}
\end{theorem}

The following result proves the convexity of the Shannon entropy power along the conjugate heat equation introduced by Perelman \cite{P1} on the Ricci flow. It can be used to characterize the shrinking Ricci solitons.

\begin{theorem}\label{thm4}  Let $(M, g(t), t\in [0, T])$ be a compact Ricci flow $\partial_t g=-2Ric_{g(t)}$. Let $u(t)$ be the fundamental solution to the conjugate heat equation 
\begin{eqnarray*}
\partial_t u=-\Delta u+Ru.\label{BHE1}
\end{eqnarray*}
Let
\begin{eqnarray*}
H(u(t))=-\int_M u\log udv, \ \ \ \ \mathcal{N}(u(t))=e^{{2\over n}H(u(t))}.
\end{eqnarray*}
Then 
\begin{eqnarray*}
{d^2 \mathcal{N} \over dt^2}={2\mathcal{N} \over n}\left[{2\over n}\left(\mathcal{F}-{n\over 2\tau}\right)^2+\int_M \left|Ric+\nabla^2 f-{g\over 2\tau }\right|^2 udv\right],  \label{NIWRF}
\end{eqnarray*}
where 
\begin{eqnarray*}
\mathcal{F}=\int_M (R+|\nabla \log u|^2)udv.
\end{eqnarray*}
In particular,  the Shannon entropy power is convex on $(0, T]$, i.e., 
\begin{eqnarray*}
{d^2\over dt^2} \mathcal{N}(u(\tau))\geq 0.
\end{eqnarray*}
Moreover,  ${d^2\over dt^2}\mathcal{N}(u(\tau))=0$ holds at some $\tau=\tau_0\in (0, T]$  if and only if  $(M, g(\tau), u(\tau))$  is a shrinking Ricci soliton 
\begin{eqnarray*}
Ric+\nabla^2 \log u={g\over 2\tau}.
\end{eqnarray*} 
\end{theorem}

As application of the concavity of the Shannon entropy power on Riemannian manifolds, we have the entropy  isoperimetric inequality  on complete Riemannian manifolds with non-negative Ricci curvature and with  maximal volume growth condition.

\begin{theorem} \label{thm5a}  Let $M$ be an $n$-dimensional complete Riemannian manifold with $Ric\geq 0$. Suppose that there exists a constant $C_n>0$ such that the maximal volume growth condition holds
$$Vol(B(x, r))\geq C_n r^n, \ \ \ \forall x\in M, r>0.$$ 
Let $u$ be the fundamental solution to the heat equation
$$
\partial_t u=\Delta u.$$
Then the following isoperimetric inequality for Shannon entropy power holds: for any probability distribution $fdv$ on $M$ such that $I(f)$ and $H(f)$ well-defined,  we have
\begin{eqnarray}
I(f)N(f)\geq \gamma_n:=4\pi e n \kappa^{2\over n},  \label{EIS1}
\end{eqnarray}
where
\begin{eqnarray*}
\kappa:=\lim\limits_{r\rightarrow \infty}\inf\limits_{x\in M} {V(B(x, r))\over \omega_n r^n},
\end{eqnarray*}
and $\omega_n$ denotes the volume of the unit ball in $\mathbb{R}^n$. Equivalently, the Stam type  logarithmic Sobolev inequality holds: for any smooth function $f$ such that $\int_M f^2dv=1$ and 
$\int_M |\nabla f|^2 dv<\infty$, we have
\begin{eqnarray*}
\int_M f^2\log f^2 dv\leq {n\over 2}\log \left({4\over \gamma_n} \int_M |\nabla f|^2dv\right). 
\end{eqnarray*}

\end{theorem}

%
%
%
%
%

\begin{remark} {After we finished the earlier version of our paper in September 2019, we found that the Entropy Differential Inequality ${\rm EDI}(K, m)$  in Theorem \ref{thm2} has been already proved by D. Bakry \cite{Ba94} on compact Riemannian manifolds, but he did not introduce the notion of the Shannon Entropy Power on manifolds.  The idea of using the Shannon entropy power to characterize the rigidity models as stated  in Theorem \ref{thm1}, Theorem \ref{thm2} and Theorem \ref{thm3} is new. Theorem \ref{thm4} uses the Shannon entropy power to characterize the shrinking Ricci solitons. It can be regarded as a natural correspondence of using Perelman's $W$-entropy to characterize the shrinking Ricci solitons. The entropy isoperimetric inequality can be also extended to complete Riemannian manifolds with $CD(0, m)$ condition and with the maximal  volume growth condition. See Section 8. 
}
\end{remark} 

The rest part of this paper is organized as follows. In Section 3, we prove the Entropy Differential Inequality (EDI) and the Entropy Power Concavity Inequality (EPCI) in Theorem \ref{thm1} and Theorem \ref{thm2} on complete Riemannian manifolds. In Section 4, we prove EDI and EPCI  in Theorem \ref{thm3} on compact $(K, m)$-super Ricci flows.  In Section 5, we prove an interesting formula between the Shannon entropy power $N$, the Fisher information $I$ and the $W$-entropy for the heat equation associated with the Laplacian or the Witten Laplacian on complete Riemannian manifolds or compact super Ricci flows. This yields an explicit  formula for the second derivative of the Shannon entropy power. In Section 6, we prove the rigidity parts in Theorem \ref{thm1}, Theorem \ref{thm2} and Theorem \ref{thm3}. In Section 7, we introduce the Shannon entropy power for the conjugate heat equation on Ricci flow and prove Theorem \ref{thm4}.  
In Section 7, we prove the entropy isoperimetric inequality in Theorem \ref{thm5a} and extend it to complete Riemannian manifolds with the $CD(0, m)$-condition and maximal volume growth condition.

\section{Concavity of the entropy power on manifolds}

In this section, we prove the concavity  of the Shannon entropy power for the heat kernel distribution associated with the Witten Laplacian on complete Riemannian manifolds with the $CD(K, m)$-condition. We need the following entropy dissipation formulas  for the heat equation associated with  the Witten Laplacian on complete Riemannian manifolds with bounded geometry condition. In the case of  compact Riemannian manifolds, it is a well-known result due to Bakry and Emery \cite{BE}.

\begin{theorem}\label{AAA}(\cite{Li12, Li16, LL15})
Let~$(M,g)$~be a complete Riemannian manifold with bounded geometry condition, and $\phi\in C^4(M)$ with $\nabla\phi\in C_b^3(M)$. Let $u$ be the fundamental solution to the heat equation ~$\partial_tu=L u$.
Let
$$H(u(t))=-\int_M u\log ud\mu.$$ Then
\begin{eqnarray}
{d\over dt}H(u(t))&=&\int_M |\nabla \log  u|^2 ud\mu, \label{ENTH1}\\
{d^2\over dt^2} H(u(t))&=&-2\int_M \Gamma_2(\nabla \log u, \nabla \log u) ud\mu, \label{ENTH2}
\end{eqnarray}
where
\begin{eqnarray*}
\Gamma_2(\nabla \log u, \nabla \log u)=\|\nabla^2\log u\|_{\rm HS}^2+Ric(L)(\nabla \log u, \nabla \log u).
\end{eqnarray*}
\end{theorem}

Moreover, we need the  Li-Yau differential Harnack inequality or the Li-Yau-Hamilton differential Harnack inequality for the heat equation associated with the Witten Laplacian on complete Riemannian manifolds with $CD(-K, m)$-condition.
When $m=n$, $\phi=0$, $L=\Delta$ and $Ric\geq -Kg$, they are due to Li-Yau \cite{LY} and Hamilton \cite{Ha}.

\begin{theorem}\label{HLYH} (\cite{LY, Ha, Li05, LL18a}) Let $(M, g)$ be a complete Riemannian manifold,  $\phi\in C^2(M)$. Suppose that there exist some
constants $m\in [n, \infty)$ and $K\geq 0$ such that
$Ric_{m, n}(L)\geq -K$. Let $u$ be a
positive solution of the heat equation $\partial_t u=Lu$.
Then the Li-Yau differential Harnack inequality holds: for all
$\alpha>1$,
\begin{eqnarray}
{|\nabla u|^2\over u^2}-\alpha {\partial_t u\over u}\leq
{m\alpha^2\over 2t}+{m\alpha^2K\over 2(\alpha-1)},\label{LYK}
\end{eqnarray}
and the Li-Yau-Hamilton differential Harnack inequality holds
\begin{eqnarray}
{|\nabla u|^2\over u^2}-e^{2Kt}{\partial_t u\over u}\leq {m\over
2t}e^{4Kt}.   \label{HmK}
\end{eqnarray}
In particular, if $Ric_{m, n}(L)\geq 0$, then the Li-Yau differential Harnack
inequality holds
\begin{eqnarray*}
{|\nabla u|^2\over u^2}-{\partial_t u\over u}\leq {m\over
2t}.
\end{eqnarray*}
\end{theorem}

\medskip

\noindent{\bf Proof of EDI$(K, m)$ and EPCI$(K, m)$ in  Theorem \ref{thm2} and Theorem \ref{thm1}}. By Theorem \ref{AAA}, we have
\begin{eqnarray*}
H'(u(t))=\int_M {|\nabla u|^2\over u}d\mu. 
\end{eqnarray*}
By \cite{Li05}, under the condition $Ric_{m, n}(L)\geq Kg$, the solution to the heat equation $\partial_t u=Lu$ is unique in $L^\infty$, and $\int_M \partial_t u d\mu=\int_M L u d\mu=0$. Thus, for any function $\alpha: [0, \infty)\rightarrow \mathbb{R}$,  we have
\begin{eqnarray*}
H'(u(t))=\int_M \left[{|\nabla u|^2\over u^2}-\alpha(t) {\partial_t u\over u}\right] ud\mu.
\end{eqnarray*}
Let $K^{-}=\max\limits\{0, -K\}$.
By the  Li-Yau Harnack estimate $(\ref{LYK})$, for any $\alpha>1$ and $t>0$, we have
\begin{eqnarray*}
H'(u(t))\leq {m\alpha^2\over 2t}+{m\alpha^2K^{-}\over 2(\alpha-1)},
\end{eqnarray*}
or,  by the Hamilton type Harnack estimate $(\ref{HmK})$, for any  $t>0$, we have 
\begin{eqnarray*}
H'(u(t))\leq {m\over 2t}e^{4K^{-}t}.
\end{eqnarray*}
From each of the above differential Harnack inequalities,  we can derive that 
\begin{eqnarray*}
\lim\limits_{t\rightarrow 0+}tH'(u(t))\leq {m\over 2}.
\end{eqnarray*}

Now we prove the entropy differential inequality.  By Theorem \ref{AAA}, we have
\begin{eqnarray*}
H'(u)=\int_M |\nabla \log u|^2 ud\mu,
\end{eqnarray*}
and
\begin{eqnarray*}
-{1\over 2}H''(u)=\int_M [\|\nabla^2\log u\|_{\rm HS}^2+Ric(L)(\nabla u, \nabla u)]ud\mu.
\end{eqnarray*}
By \cite{BE, Li05}, it holds
\begin{eqnarray*}
\|\nabla^2\log u\|_{\rm HS}^2+Ric(L)(\nabla u, \nabla u)\geq {|L\log u|^2\over m}+Ric_{m, n}(L)(\nabla \log u, \nabla \log u).
\end{eqnarray*}
Integrating the above inequality on $M$ with respect to $u\mu$, and using the second entropy dissipation formula in Theorem \ref{AAA}, we get
\begin{eqnarray*}
-{1\over 2}H''(u)\geq  \int_M\left[{|L\log u|^2\over m}+Ric_{m, n}(L)(\nabla \log u, \nabla\log u)\right]ud\mu.
\end{eqnarray*}
By the Cauchy-Schwarz inequality and integration by parts,  we have
\begin{eqnarray*}
\int_M |L\log u|^2 ud\geq \left(\int_M L\log u ud\mu\right)^2=\left(\int_M |\nabla \log u|^2 ud\mu\right)^2.
\end{eqnarray*}
This proves the entropy differential inequality $(\ref{EDI})$.

Finally, note that
\begin{eqnarray*}
N'&=&{2H'\over m}N,\\
N''&=&{2H''\over m}N+{4H'^2\over m^2}N={2N\over m}(H''+{2H'^2\over m})\leq -2K H' {2N\over m}.
\end{eqnarray*}
Thus
\begin{eqnarray*}
N''\leq-2KN'.
\end{eqnarray*}
The proof of  EDI$(K, m)$ and EPCI$(K, m)$ in Theorem \ref{thm2} is completed. Taking $m=n$ and $\phi=0$ in Theorem \ref{thm2}, we can derive  EDI$(K, n)$ and EPCI$(K, n)$ in Theorem \ref{thm1}. \hfill $\square$

\section{Concavity of the entropy power on super Ricci flows}

In this section, we prove the concavity  of the Shannon entropy power for the heat distribution associated with the Witten Laplacian on compact $(K, m)$-super Ricci flows. We need the following entropy dissipation formulas  for the heat equation associated with the time dependent Witten Laplacian on complete Riemannian manifolds with super Ricci flows. 

\begin{theorem}\label{BBB}(\cite{Li12, Li16, LL15})
Let~$(M,g)$~be a compact Riemannian manifold with a family of time dependent metrics $g(t)$ and potentials $\phi(t)\in C^4(M)$, $t\in [0, T]$. Let $u$ be a positive solution to the heat equation ~$\partial_tu=L u$.
Let
$$H(u(t))=-\int_M u\log ud\mu.$$ Then
\begin{eqnarray}
{d\over dt}H(u(t))&=&\int_M |\nabla \log  u|^2 ud\mu, \label{ENTH1}\\
{d^2\over dt^2} H(u(t))&=&-2\int_M \Gamma_2(\nabla \log u, \nabla \log u) ud\mu, \label{ENTH2}
\end{eqnarray}
where
\begin{eqnarray*}
\Gamma_2(\nabla \log u, \nabla \log u)=\|\nabla^2\log u\|_{\rm HS}^2+\left({1\over 2}{\partial g\over \partial t}+Ric(L)\right)(\nabla \log u, \nabla \log u).
\end{eqnarray*}
\end{theorem}

\medskip

\noindent{\bf Proof of  EDI$(K, m)$ and EPCI$(K, m)$ in Theorem \ref{thm3}}. By the same argument as used in the proof of Theorem \ref{thm1},  we can prove that 
\begin{eqnarray*}
\Gamma_2(\nabla \log u, \nabla \log u)\geq {|L\log u|^2\over m}+\left({1\over 2}{\partial g\over \partial t}+Ric_{m, n}(L)\right)(\nabla \log u, \nabla \log u).
\end{eqnarray*}
Integrating the above inequality on $M$ with respect to $u\mu$, and using Theorem \ref{BBB}, we get
\begin{eqnarray*}
-{1\over 2}H''(u)\geq  \int_M\left[{|L\log u|^2\over m}+\left({1\over 2}{\partial g\over \partial t}+Ric_{m, n}(L)\right)(\nabla \log u, \nabla\log u)\right]ud\mu.
\end{eqnarray*}
The rest of the proof of   ECI$(K, m)$ and EPCI$(K, m)$  in Theorem \ref{thm3}  is similar to the one of Theorem \ref{thm2}.   \hfill $\square$

\medskip

\section{Shannon entropy power and $W$-entropy}

In this section, we prove an interesting formula between the Shannon entropy power, the Fisher information and the $W$-entropy for the heat equation associated with the Laplacian or  the Witten Laplacian on manifolds and on $(K, m)$-super Ricci flows. This derives an explicit formula  for the second derivative of the Shannon entropy power which enables us to prove the rigidity theorem of the Shannon entropy power  in Section 6.

 Let 
\begin{eqnarray}
H_m(u(t))=H(u(t))-{m\over 2}\log(4\pi et). \label{Hm}
\end{eqnarray}
Inspired by Perelman \cite{P1}, we introduce the $W$-entropy by the Boltzmann entropy formula (see \cite{Li12, Li16, LL15})
\begin{eqnarray}
W_m(u(t))={d\over dt} (tH_m(u(t)).\label{Wm}
\end{eqnarray}
By  \cite{Li12, Li16, LL15},  we have 
\begin{eqnarray*}
W_m(u(t))=\int_M \left(t |\nabla f|^2+f-m\right)u d\mu.
\end{eqnarray*}
Moreover, the following $W$-entropy formula has been proved in  \cite{Li12, Li16, LL15}
\begin{eqnarray} 
{d \over dt}W_m(u(t))&=&-2\int_M
\left(t \left|\nabla^2 f-{g\over 2t}\right|^2+Ric_{m, n}(L)(\nabla f, \nabla f)\right)ud\mu \nonumber\\
& & \hskip2cm -{2\over m-n}\int_M t \left({\nabla\phi\cdot\nabla f}+{m-n\over 2t}\right)^2ud\mu.  \label{derW}
\end{eqnarray}
In the case $m=n$, $\phi=0$ and $L=\Delta$, this formula was due to  Ni \cite{N1}.  

Note that
\begin{eqnarray}
W_m(u(t))=H(u(t))-{m\over 2}\log(4\pi e t)+tH'(u(t))-{m\over 2},
\end{eqnarray}
and
\begin{eqnarray}
{d\over dt}W_m(u(t))=tH''(u(t))+2H'(u(t))-{m\over 2t}.
\end{eqnarray}
Hence

%
%
\begin{eqnarray*}
H''+{2\over m}H'^2
={2\over m}\left(H'-{m\over 2t}\right)^2+{1\over t}{d\over dt}W_m(u(t)).
\end{eqnarray*}   
By the fact that
\begin{eqnarray*}
N''={2N\over m}\left(H''+{2\over m}H'^2\right),
\end{eqnarray*}
we prove the following NIW formula which has its own interest.

\begin{theorem}\label{ThNIW}  The following NIW formula  holds

 \begin{eqnarray}
{d^2 N\over dt^2}={2N\over m}\left[{2\over m}\left(I-{m\over 2t}\right)^2+{1\over t}{dW_m\over dt}\right].  \label{NIW}
\end{eqnarray}
Moreover, 
\begin{eqnarray}
{m\over 2N}{d^2 N\over dt^2}&=&-2\int_M Ric_{m, n}(L)(\nabla f, \nabla f)ud\mu-{2\over m}\int_M \left[ Lf-\int_M Lf ud\mu\right]^2u d\mu \nonumber\\
& &-2\left({1\over n}-{1\over m}\right)\int_M \left[\Delta f+{n\over m-n}\nabla \phi\cdot\nabla f\right]^2 ud\mu-2\int_M \left\|\nabla^2 f-{\Delta f\over n}g\right\|_{\rm HS}^2ud\mu.\nonumber\\
& & \label{NN1}
\end{eqnarray}
In particular, under the condition $Ric_{m, n}(L)\geq Kg$, we have
\begin{eqnarray*}
{d^2 N\over dt^2}\leq-2K{dN\over dt}.
\end{eqnarray*}
Equivalently, the entropy power concavity inequality EPCI$(K, m)$ holds.
\end{theorem}
{\it Proof}. It remains to prove $(\ref{NN1})$. Combining $(\ref{NIW})$  with the $W$-entropy formula $(\ref{derW})$, we have
\begin{eqnarray}
{m\over 2N}{d^2 N\over dt^2}&=&{2\over m}\left(I-{m\over 2t}\right)^2-2\int_M \left(\left\|\nabla^2 f-{g\over 2t}\right\|_{\rm HS}^2+Ric_{m, n}(L)(\nabla f, \nabla f)\right)ud\mu \nonumber\\
& & \hskip2cm -{2\over m-n}\int_M  \left({\nabla\phi\cdot\nabla f}+{m-n\over 2t}\right)^2ud\mu. \label{NNN1}
\end{eqnarray}
Using 
\begin{eqnarray*}
\|A\|_{\rm HS}^2={|{\rm Tr A}|^2\over n}+\left\|A-{{\rm Tr}A\over n} g\right\|^2_{\rm HS}
\end{eqnarray*}
to  symmetric matrix $A=\nabla^2 f-{g\over 2t}$, and applying 
\begin{eqnarray*}
(a+b)^2 = {a^2\over 1+\varepsilon}-{b^2\over \varepsilon}+{\varepsilon\over 1+\varepsilon}\left(a+{1+\varepsilon\over \varepsilon}b\right)^2
\end{eqnarray*}
to $a=Lf-{m\over 2t}$ and $b=\nabla \phi\cdot\nabla f+{m-n\over 2t}$ with $\varepsilon={m-n\over n}$, we have
\begin{eqnarray*}
\left\|\nabla^2 f-{g\over 2t}\right\|_{\rm HS}^2
&=&{1\over n} \left|\Delta f-{n\over 2t}\right|^2+\left\|\nabla^2 f-{\Delta f\over n}g\right\|_{\rm HS}^2\\
&=&{1\over n} \left|Lf-{m\over 2t}+\left(\nabla \phi\cdot\nabla f+{m-n\over 2t}\right)   \right|^2+\left\|\nabla^2 f-{\Delta f\over n}g\right\|_{\rm HS}^2\\\
&=& {1\over m} \left|Lf-{m\over 2t}  \right|^2-{1\over m-n}\left(\nabla \phi\cdot\nabla f+{m-n\over 2t}\right)^2\\
& &+{m-n\over mn}\left[Lf+{m\over m-n}\nabla \phi\cdot\nabla f\right]^2+\left\|\nabla^2 f-{\Delta f\over n}g\right\|_{\rm HS}^2.
\end{eqnarray*}
Note that
\begin{eqnarray*}
I=\int_M |\nabla \log u|^2 ud\mu=\int_M Lf ud\mu,
\end{eqnarray*}
and
\begin{eqnarray*}
\left(I-{m\over 2t}\right)^2-\int_M \left|Lf-{m\over 2t}\right|^2 ud\mu
=-\int_M \left[ Lf-\int_M Lf ud\mu\right]^2u d\mu.
\end{eqnarray*}
Combining these with $(\ref{NNN1})$, we derive $(\ref{NN1})$.  The EPCI$(K, m)$ follows. \hfill $\square$

\medskip

In particular, when $m=n$ and $\phi=C$, we have the following

\begin{theorem}\label{NIW2}  The following NIW formula  holds
 \begin{eqnarray}
{d^2 N\over dt^2}={2N\over n}\left[{2\over m}\left(I-{n\over 2t}\right)^2+{1\over t}{dW_n\over dt}\right].  \label{NIWn}
\end{eqnarray}
Moreover, 
\begin{eqnarray}
{n\over 2N}{d^2 N\over dt^2}
&=&-2\int_M Ric(\nabla f, \nabla f) ud\mu-{2\over n}\int_M\left[\Delta f-\int_M (\Delta f) udv\right]^2 udv\nonumber\\
& &\hskip2cm -2\int_M \left\|\nabla^2 f-{\Delta f\over n}g\right\|_{\rm HS}^2 udv.\label{NIWn2}
\end{eqnarray}
In particular, under the condition $Ric\geq Kg$, we have
\begin{eqnarray*}
{d^2 N\over dt^2}\leq-2K{dN\over dt}.
\end{eqnarray*}
Equivalently, the entropy power concavity inequality EPCI$(K, n)$ holds.
\end{theorem}

\medskip

Similarly, we can prove the following result which extends Theorem \ref{ThNIW} to (weighted) complete Riemannian manifolds with  Ricci curvature (or  the $m$-dimensional Bakry-Emery Ricci curvature) bounded from below by a constant.  To save the length of the paper, we omit the proof. 

\begin{theorem}\label{NIWT} The following formula holds for the Fisher information $I=H'$, Shannon entropy power $N$ and Perelman's $W$-entropy for the heat equation $\partial_t u=Lu$ on complete Riemannian manifolds with $CD(K, m)$-condition or compact $(K, m)$-super Ricci  flows
\begin{eqnarray}
{d^2 N\over dt^2}+2K{d N\over dt}={2N\over m}\left[{2\over m}\left(I-{m(1+Kt)\over 2t}\right)^2+{1\over t}{dW_{m, K}\over dt}\right].  \label{NIWK}
\end{eqnarray}
where 
\begin{eqnarray}
W_{m, K}=\int_M t\left(|\nabla f|^2+f-m\left(1+{Kt\over 2}\right)^2\right)ud\mu.
\end{eqnarray}
Moreover, we have
\begin{eqnarray}
& &{m\over 2N}\left[{d^2 N\over dt^2}+2K{dN\over dt}\right]\nonumber\\
&=&-{2\over m}\int_M \left[ Lf-\int_M Lf ud\mu\right]^2u d\mu-2\int_M \left({1\over 2}{\partial g\over \partial t}+Ric_{m, n}(L)-Kg\right)(\nabla f, \nabla f)ud\mu \nonumber\\
& &-2\left({1\over n}-{1\over m}\right)\int_M \left[\Delta f+{n\over m-n}\nabla \phi\cdot\nabla f\right]^2 ud\mu-2\int_M \left\|\nabla^2 f-{\Delta f\over n}g\right\|_{\rm HS}^2ud\mu.\nonumber\\
& & \label{NIWt}
\end{eqnarray}
In particular, under the condition $Ric_{m, n}(L)\geq Kg$ (in time independent case)  or compact $(K, m)$-super Ricci flow (time dependent case), i.e., 
\begin{eqnarray}
{1\over 2}{\partial g\over \partial t}+Ric_{m, n}(L)\geq Kg, 
\end{eqnarray} we have
\begin{eqnarray*}
{d^2 N\over dt^2}\leq-2K{dN\over dt}.
\end{eqnarray*}
Equivalently, the entropy power concavity inequality EPCI$(K, m)$ (i.e., $(\ref{ECI1})$) holds.
\end{theorem}
\section{Rigidity theorems for the Shannon entropy power}

Now we prove the rigidity theorems for the Shannon entropy power on complete Riemannian manifolds with $CD(K, m)$-conditions  and on compact $(K, m)$-super Ricci flows. 

\medskip

\noindent{\bf Proof of rigidity part in Theorem \ref{thm1}}. 
By $(\ref{NIWn2})$ in Theorem \ref{NIW2}, in the case $Ric\geq Kg$,  assuming that ${d^2N\over dt^2}=-2K{dN\over dt}$ holds on $(0, T]$ for some $T>0$, then 
\begin{eqnarray*}
Ric(\nabla f, \nabla f)=K|\nabla f|^2, \ \ \ \Delta f=\int_M (\Delta f)udv, \ \ \ \nabla^2 f={\Delta f\over n}g.
\end{eqnarray*}
By the Varadhan type small time asymptotic behavior of the heat kernel of $\partial_t =\Delta u$ on complete Riemannian manifolds
\begin{eqnarray*}
\lim\limits_{t\rightarrow 0+} t\log p_t(o, x)=-{d^2(o, x)\over 4}, \ \ \ \ \forall x\in M,
\end{eqnarray*}
where $o\in M$ is a fixed referenced point. This yields
\begin{eqnarray*}
tf(x, t)= {d^2(o, x)\over 4}+o(1),   \ \ \ \   t\rightarrow 0+.
\end{eqnarray*}
This together with the first identity implies that $Ric=Kg$ as $t\nabla f(x, t)$ spans $T_xM$.  The second one implies that implies that $\Delta f=\int_M |\nabla \log u|^2 udv=H_K(t)>0$ is a positive constant depending only on $t$ and the reference point $o$. The third identity implies that \begin{eqnarray*}
\nabla^2 f(x, t)={H'_K(t)\over n}g.
\end{eqnarray*}
When $K=0$, by the same argument as used in \cite{Li12},  see also the proof of rigidity part in Theorem \ref{thm2} below,  we can prove that $M$ is isometric to $\mathbb{R}^n$, and 
$$f(x, t)={\|x\|^2\over 4t}, \ \ \ \forall~ x\in \mathbb{R}^n, t>0.$$
The proof of Theorem \ref{thm1}  is completed. \hfill $\square$

\medskip

\noindent{\bf Proof of  rigidity part in Theorem \ref{thm2} and Theorem \ref{thm3}}. 
By $(\ref{NN1})$ in Theorem \ref{ThNIW},  ${d^2N\over dt^2}=0$ holds on $(0, T]$ for some $T>0$ if and only if 
\begin{eqnarray*}
& &Ric_{m, n}(L)(\nabla f, \nabla f)=K|\nabla f|^2, \ \ \  L f=\int_M (L f)ud\mu,\\
& & \ \ \ \  \Delta f+{n\over m-n}\nabla f\cdot \nabla \phi=0, \ \ \ \ \ \ \ \ \nabla^2 f={\Delta f\over n}g.
\end{eqnarray*}
Similarly to the proof of the rigidity part of Theorem \ref{thm1}, by the Varadhan type small time asymptotic behavior of the heat kernel of $\partial_t u=Lu$ on complete Riemannian manifolds with the $CD(K, m)$-condition, which was proved by the second author in \cite{Li12}, 
 the first identity implies that $Ric_{m, n}(L)=Kg$. Moreover, we can prove that  $Lf=\int_M |\nabla \log u|^2 ud\mu=H'_K(t)$, $\Delta f={n\over m}H'_K(t)$ and $\nabla^2 f={H'_K(t)\over m}g$.  

In particular, when $K=0$, by the same argument as above, $L f=\int_M |\nabla \log u|^2 u d\mu=H'(t)>0$,  $\Delta f={n\over m}H'(t)$, and 
\begin{eqnarray*}
\nabla^2 f={\Delta f\over n} g={H'(t)\over m}g.
\end{eqnarray*}
Thus, $f$ is a strict convex function on $M$, and $M$ must be diffeomorphic to $\mathbb{R}^n$.   Fix $t=t_0$ and let $x_0$ be the minimal point of $f(
\cdot, t)$. Integrating along the shortest geodesic linking $x_0\in M$ and any $x\in M$, we have
\begin{eqnarray*}
f(x, t_0)-f(x_0, t_0)={H'(t_0)\over 2m} r^2(x_0, x). 
\end{eqnarray*}
This implies  
\begin{eqnarray*}
\Delta r^2(x, x_0)={2m\over H'(t_0)}\Delta f(x, t_0)=2n.
\end{eqnarray*}
Hence,  $M$ is isometric to $\mathbb{R}^n$.  It follows from $Ric_{m, n}(L)=Ric+\nabla^2\phi-{\nabla\phi\otimes \nabla\phi\over m -n}=0$ that 
\begin{eqnarray*}
\nabla^2\phi-{\nabla\phi\otimes \nabla\phi\over m -n}=0. 
\end{eqnarray*}
Taking trace on the both sides of the above equality, we have 
\begin{eqnarray*}
\Delta\phi-{|\nabla\phi|^2\over m -n}=0. 
\end{eqnarray*}
Thus
\begin{eqnarray*}
\Delta e^{-{\phi\over m-n}}=0.
\end{eqnarray*}
By the classical Strong Liouville theorem on $\mathbb{R}^n$, we can conclude that $e^{-{\phi\over m-n}}$ must be a constant, which yields that  $\phi$ must be a constant. Indeed, by the generalized Cheeger–Gromoll splitting theorem proved by the second author with Fang and Zhang (see Theorem 1.3, p. 565 in \cite{FLZ}), we can also 
derive that $\phi$ must be a constant.  It follows that 
\begin{eqnarray*}
u(x,t) = {e^{-{\|x\|^2\over 4t}}\over (4\pi t)^{n/2}}, \  \ \ \forall  t>0, x\in M=\mathbb{R}^n,
\end{eqnarray*}
and  one can chose $m=n$. This yields that  $f(x, t)={\|x\|^2\over 4t}$ for all  $t>0$, $x\in M=\mathbb{R}^n$ and $H(t)={n\over 2}\log(4\pi et)$ and $H'(t)={n\over 2t}$. 
This finishes the proof of Theorem \ref{thm4}. 
\medskip

Similarly, the rigidity part in Theorem \ref{thm3} can be proved by $(\ref{NIWt})$ in Theorem \ref{NIWT}.   \hfill $\square$ 

\section{Shannon entropy power on Ricci flow}

In this section we prove the convexity of the Shannon entropy power along the conjugate heat equation introduced by Perelman \cite{P1} on the Ricci flow. It can be used to characterize the shrinking Ricci solitons.

Let $M$ be a compact manifold equipped with the Ricci flow 
\begin{eqnarray}
\partial_t g(t)=-2Ric_{g(t)}.\label{RF}
\end{eqnarray}
Let $u(t)={e^{-f}\over (4\pi t)^{n/2}}$ be the fundamental solution to the backward heat equation 
\begin{eqnarray}
\partial_t u=-\Delta u+Ru.\label{BHE}
\end{eqnarray}
Let $\tau$ be such that 
$$\partial_t \tau=-1.$$

Let 
\begin{eqnarray}
\mathcal{H}(g, u)=-\int_M u\log u dv
\end{eqnarray}
be the Shannon entropy.  
In \cite{P1}, Perelman introduced the following $\mathcal{F}$-function and $\mathcal{W}$-entropy
\begin{eqnarray}
\mathcal{F}(g, u)&=&\int_M (R+|\nabla \log u|^2)udv,\\
\mathcal{W}(g, u))&=&\int_M \left(\tau (R+|\nabla f|^2)+f-n\right)u dv,
\end{eqnarray} 
and proved the following beautiful  formulas 
\begin{eqnarray} 
{d\over dt} \mathcal{H}(g(t), u(t))&=&-\mathcal{F}(g(t), u(t)), \label{FH}\\
{d\over dt} \mathcal{F}(g(t), u(t))&=&2\int_M |Ric+\nabla^2 \log u|^2 udv, \label{Ft}\\
{d \over dt}\mathcal{W}(g(\tau), (u(\tau))&=&2\tau \int_M \left|Ric+\nabla^2 f-{g\over 2\tau }\right|^2 udv. \label{Wt}
\end{eqnarray}
  
In \cite{Li12},  a probabilistic interpretation was given for the $W$-entropy on Ricci flow. Indeed, let 
\begin{eqnarray}
\mathcal{H}_n(g(\tau), u(\tau))=\mathcal{H}(g(\tau), u(\tau))-{n\over 2}\log(4\pi e\tau) \label{Hn}
\end{eqnarray}
be the difference of the Shannon entropy $\mathcal{H}(g(\tau), u(\tau))$ of the heat kernel measure $u(t)dv_{g(t)}$  on Ricci flow and the Shannon entropy $H(\overline{u})={n\over 2} \log(4\pi et)$ of the Gaussian heat kernel 
$\overline{u}(t)={1\over (4\pi t)^{n/2}}
e^{-{\|x\|^2\over 4t}}$ on $\mathbb{R}^n$. 
Then Perelman's $W$-entropy for the Ricci flow satisfies  the Boltzmann entropy formula in statistical mechanics
\begin{eqnarray}
\mathcal{W}(g(\tau), u(\tau))=-{d\over d\tau} (\tau \mathcal{H}_n(g(\tau), u(\tau)).\label{W}
\end{eqnarray}

Note that
\begin{eqnarray}
\mathcal{W}=-\mathcal{H}+{n\over 2}\log(4\pi e \tau)-\tau{d\mathcal{H}\over d\tau} +{n\over 2},
\end{eqnarray}
and
\begin{eqnarray}
{d\over d\tau}\mathcal{W}=-\tau{d^2 \mathcal{H}\over d\tau^2}-2{d\mathcal{H}\over d\tau}+{n\over 2\tau}.
\end{eqnarray}
Hence

\begin{eqnarray*}
{d^2 \mathcal{H}\over d\tau^2}=-{1\over \tau}{d\over d\tau}\mathcal{W}-{2\over \tau}{d\mathcal{H}\over d\tau}+{n\over 2\tau^2},
\end{eqnarray*}
which yields 
\begin{eqnarray*}
\mathcal{H}''+{2\over n}\mathcal{H}'^2={2\over n}\left(\mathcal{H}'-{n\over 2\tau}\right)^2-{1\over \tau}{d\over d\tau}\mathcal{W}(u(t)).
\end{eqnarray*}
where $\mathcal{H}'={d\mathcal{H}\over d\tau}$, and ${\mathcal H}''={d^2\mathcal{H}\over d\tau^2}$.  

The Shannon entropy power for the conjugate heat equation is defined as follows
\begin{eqnarray}
\mathcal{N}=\mathcal{N}(g(\tau), u(\tau))=e^{{2\over n}\mathcal{H}(g(\tau), u(\tau))}.
\end{eqnarray}
Obviously, we have
\begin{eqnarray*}
{d^2\mathcal{N}\over d\tau^2}={2\mathcal{N}\over n}\left(\mathcal{H}''+{2\over n}\mathcal{H}'^2\right).
\end{eqnarray*}
Note that
\begin{eqnarray*}
{d\over d\tau}=-{d\over dt}, \ \ \ \  {d^2\over d\tau^2}={d^2\over dt^2}.
\end{eqnarray*}
This proves the NFW formula in Theorem \ref{thm4}.  More precisely, 
\begin{eqnarray*}
{d^2\mathcal{N}\over d\tau^2}={2\mathcal{N}\over n}\left[{2\over n}\left(\mathcal{F}-{n\over 2\tau}\right)^2+
{1\over \tau}{d\mathcal{W}\over dt}\right].
\end{eqnarray*}

Note that ${d\mathcal{W}\over dt}\geq 0$. Hence ${d^2\mathcal{N}\over d\tau^2}\geq 0$. This proves the convexity of the Shannon entropy power along the conjugate heat equation on the Ricci flow. 

Moreover,  if ${d^2\mathcal{N}\over d\tau^2}=0$ holds at some $\tau=\tau_0$, then we must have
\begin{eqnarray*}
\mathcal{F}={n\over 2\tau},
\end{eqnarray*}
and
\begin{eqnarray*}
{d\mathcal{W}\over dt}=0.
\end{eqnarray*}
By Perelman's $W$-entropy formula $(\ref{Wt})$, the second identity implies that $M$ is a shrinking Ricci soliton
\begin{eqnarray*}
Ric+\nabla^2f ={g\over 2\tau}.
\end{eqnarray*}
Taking trace on the both side, we have
\begin{eqnarray*}
R+\Delta f ={n\over 2\tau}.
\end{eqnarray*}
Integration by parts yields
\begin{eqnarray*}
\mathcal{F}=\int_M (R-\Delta \log u) udv=\int_M (R+\Delta f)  udv={n\over 2\tau}.
\end{eqnarray*}
The proof of Theorem \ref{thm4} is completed. \hfill $\square$

\section{Entropy Isoperimetric Inequality on manifolds}

In this section, as an application of the concavity of the Shannon entropy power for heat equation, we prove an entropy isoperimetric inequality  on 
complete Riemannian manifolds with  non-negative Ricci curvature or $m$-dimensional Bakry-Emery Ricci curvature and  maximal volume growth condition. 

First, we consider the case of complete Riemannian manifolds with non-negative Ricci curvature. In this cas, the entropy power concavity inequality EPCI$(0, n)$  (see Theorem \ref{thm1})
\begin{eqnarray*}
{d^2\over dt^2}N(u(t))\leq 0
\end{eqnarray*}
implies that ${d\over dt}N(u(t))$ is nonincreasing in $t$ and hence $\lim\limits_{t\rightarrow \infty}{d\over dt}N(u(t))$ exists. Moreover,  using the first order entropy dissipation formula in Theorem \ref{AAA}, we have

\begin{eqnarray*}
{d\over dt}N(u(t))\leq \left. {d\over dt} \right|_{t=0} N(u(t))={2\over n}I(f)N(f).  
\end{eqnarray*}
Thus 
\begin{eqnarray}
 {2\over n}I(f)N(f)\geq \lim\limits_{t\rightarrow \infty}  {d\over dt} N(u(t)). \label{DNI1}
\end{eqnarray}

Again, the first order entropy dissipation formula in Theorem \ref{AAA}  implies 

\begin{eqnarray*}
 {d\over dt}N(u(t))={2\over n} I(u(t)) N(u(t))\geq 0.
\end{eqnarray*}
Thus $N(u(t))$ is nodecreasing in $t$ and we have
\begin{eqnarray}
\lim \limits_{t\rightarrow \infty}{d\over dt}N(u(t))= \lim\limits_{t\rightarrow \infty}{N(u(t))\over t}.  \label{DNI2}
\end{eqnarray}
It remains to prove that 
$\lim\limits_{t\rightarrow \infty}{N(u(t))\over t}$ is finite and find its exact value in terms of geometric constant of manifolds. 

Let
$$
H_n(u(t))=H(u(t))-{n\over 2}\log(4\pi et).
$$ Then 
\begin{eqnarray}
{N(u(t))\over t}
=(4\pi e)e^{{2\over n}H_n(u(t))}.  \label{DNI3}
\end{eqnarray}

Under the condition $Ric\geq 0$, the Li-Yau Harnack inequality yields
\begin{eqnarray*}
{d\over dt}H_n(u(t))=-\int_M \left[\Delta \log u+{n\over 2t}\right] u d\mu\leq 0.
\end{eqnarray*}
Thus the limit $\lim\limits_{t\rightarrow \infty} H_n(u(t))$ exists and we need only to prove  $\lim\limits_{t\rightarrow \infty}H_n(u(t))$ is finite and find its exact value in terms of geometric constant of manifolds. 

The following result was proved by L. Ni \cite{N1} using sharp bound of heat kernel estimate on complete Riemannian manifolds with non-negative Ricci curvature. 

\begin{proposition} \label{Ni} (\cite{N1}) Let $M$ be a complete Riemannian manifold with $Ric\geq 0$ and the maximal volume growth condition: there exists a constant $C_n>0$ such that 
$$Vol(B(x, r))\geq C_n r^n, \ \ \ \forall x\in M, r>0.$$ 
 Let $u(t, x)=p_t(x, o)$ be the fundamental solution to the heat equation $\partial_t u=\Delta u$, where $o\in M$ is fixed.  Then 
\begin{eqnarray}
\lim\limits_{t\rightarrow \infty} H_n(u(t))=\log \kappa,  \label{DNI4}
\end{eqnarray}
where
\begin{eqnarray*}
\kappa:=\lim\limits_{r\rightarrow \infty}\inf\limits_{x\in M} {V(B(x, r))\over \omega_n r^n},
\end{eqnarray*}
and $\omega_n$ denotes the volume of the unit ball in $\mathbb{R}^n$. 
\end{proposition}

\noindent{\bf Proof of Theorem \ref{thm5a}}.  The entropy isoperimetric inequality $(\ref{EIS1})$ follows from  $(\ref{DNI1})$,  $(\ref{DNI2})$,  $(\ref{DNI3})$ and  $(\ref{DNI4})$ in Proposition \ref{Ni}. 
In general, let $\gamma_n$ be a positive constant and assume the following entropy  isoperimetric inequality  holds
\begin{eqnarray}
I(f)N(f)\geq \gamma_{n}. \label{EIIN}
\end{eqnarray}
Then, equivalently, the Stam logarithmic Sobolev inequality (LSI)  holds
\begin{eqnarray*}
\int_{M} f\log f dv\leq {n\over 2}\log \left({1\over \gamma_n} \int_{M}{|\nabla f|^2\over f}dv\right).
\end{eqnarray*}
Replacing $f$ by $f^2$ with $\int_{M} f^2dv=1$ and $\int_M |\nabla f|^2dv<\infty$,  the above LSI is equivalent to the Stam type logarithmic Sobolev inequality
\begin{eqnarray}
\int_{M} f^2\log f^2 dx\leq {n\over 2}\log \left({4\over \gamma_n} \int_M |\nabla f|^2dv\right). \label{LSI2a}
\end{eqnarray}

\hfill $\square$


\begin{remark}{
Note that when $M=\mathbb{R}^n$, $\kappa=1$, the above inequality reads 
\begin{eqnarray}
N(f)I(f)\geq 2\pi en. \label{NXY5}
\end{eqnarray}
Indeed, the  Entropy Isoperimetric Inequality $(\ref{NXY5})$  is equivalent to the following  logarithmic Sobolev inequality which was first proved by 
Stam \cite{Stam}: for any smooth probability density function $f$ on $\mathbb{R}^n$ with $\int_{\mathbb{R}^n} {|\nabla f|^2\over f} dx<\infty$, it holds
 \begin{eqnarray*}
\int_{\mathbb{R}^n} f\log f dx\leq {n\over 2}\log \left({1\over 2\pi e n} \int_{\mathbb{R}^n}{|\nabla f|^2\over f}dx\right). \label{LSI1}
\end{eqnarray*}
Replacing $f$ by $f^2$ with $\int_{\mathbb{R}^n} f^2dx=1$,  the above LSI   is equivalent to the Stam type logarithmic Sobolev inequality: for any smooth $f$ on $\mathbb{R}^n$ with $\int_{\mathbb{R}^n} f^2dx=1$ and $\int_{\mathbb{R}^n} |\nabla f|^2 dx<\infty$, it holds

\begin{eqnarray*}
\int_{\mathbb{R}^n} f^2\log f^2 dx\leq {n\over 2}\log \left({2\over \pi e n} \int_{\mathbb{R}^n} |\nabla f|^2dx\right). 
\end{eqnarray*}
We would like to mention that, if $M$ is an $n$-dimensional complete Riemannian manifold with non-negative Ricci curvature on which the logarithmic Sobolev inequality $(\ref{LSI2a})$ 
holds with the same constant $\gamma_n=2\pi e n$ as on the $n$-dimensional Euclidean space,  then $M$ must be  isometric to $\mathbb{R}^n$. 
See \cite{BCL, N1}.}
\end{remark}

The above argument can be extended to general case of weighted complete Riemnanian manifolds with the $CD(0, m)$ and maximal volume growth conditions. Indeed, by similar argument as used for the proof of Proposition \ref{Ni}, and based on two-sides heat kernel estimates and the maximal volume growth property, H. Li \cite{HLi} extended Proposition \ref{Ni} to the so-called $RCD(0, N)$ metric measure spaces with maximum volume growth condition for $N\in \mathbb{N}$ with $N\geq 2$. Thus, as it is well-known that  all 
weighted complete Riemannian manifolds with the $CD(0, m)$-condition are $RCD(0, N)$ metric measure space with $N=m\geq 2$,  we can use the entropy dissipation formula in Theorem \ref{AAA},  the entropy power concavity inequality in Theorem \ref{thm2} and the extended version of Proposition \ref{Ni} on weighted complete Riemannian manifolds with $CD(0, m)$ and maximal volume growth conditions to prove the following isoperimetric inequality for Shannon entropy power on weighted complete Riemannian manifolds.  

\begin{theorem} \label{thm6} Let $M$ be an $n$-dimensional complete  Riemannian manifolds with $Ric_{m, n}(L)\geq 0$,  where $m\in \mathbb{N}$ with $m>n$. Suppose that there exists a constant $C_m>0$ such that the maximal volume growth condition holds
\begin{eqnarray*}
\mu(B(x, r))\geq C_{m} r^m, \ \ \ \forall x\in M, r>0. \label{mvg}
\end{eqnarray*}
Let $u$ be the fundamental solution to the heat equation
$$
\partial_t u=L u.$$
Then
\begin{eqnarray*}
\lim\limits_{t\rightarrow \infty} H_m(u(t))=\log \kappa,
\end{eqnarray*}
where
\begin{eqnarray*}
\kappa:=\lim\limits_{r\rightarrow \infty}\inf\limits_{x\in M} {V(B(x, r))\over \omega_m r^m},
\end{eqnarray*}
and $\omega_m$ denotes the volume of the unit ball in $\mathbb{R}^m$. 
Furthermore, the following entropy  isoperimetric inequality holds: for any probability distribution $fd\mu$ such that $I(f)$ and $H(f)$ are well-defined, we have
\begin{eqnarray*}
I(f)N(f)\geq \gamma_m:=2\pi e m\kappa^{2\over m}.
\end{eqnarray*}Equivalently, the Stam type  logarithmic Sobolev inequality holds: for any smooth  function $f$ such that $\int_M f^2d\mu=1$ and $\int_M |\nabla f|^2 d\mu<\infty$, we have
\begin{eqnarray*}
\int_M f^2\log f^2 d\mu\leq {m\over 2}\log \left({4\over \gamma_m} \int_M |\nabla f|^2d\mu\right). \label{LSI3}
\end{eqnarray*}

\end{theorem}
%
%

In general case of complete Riemannian manifolds with the $CD(K, m)$-condition or $(K, m)$-super Ricci flows, we have the following entropy isoperimetric inequality. 

\begin{theorem} Let $M$ be an $n$-dimensional complete Riemannian manifold with $Ric_{m, n}(L)\geq Kg$ for some constants $m\geq n$ and $K\in \mathbb{R}$, and $Q(u)=N(u)I(u)$ be the entropy isoperimetric profile  along the heat equation $\partial_t u=L u$. Then 
\begin{eqnarray*}
{d\over dt} Q(u(t))\leq -2K Q(u(t)).
\end{eqnarray*}
For all $0\leq s<t<\infty$, we have
\begin{eqnarray*}
e^{2Kt} Q(u(t))\leq e^{2Ks}Q(u(s)).  \label{Qust}
\end{eqnarray*}
In particular, when $K>0$, we have
\begin{eqnarray*}
Q(u(t))\leq e^{-2Kt}Q(u(0)).
\end{eqnarray*} 
The same conclusion holds on compact $(K, m)$-super Ricci flows  as in Theorem \ref{thm3}.
\end{theorem}
{\it Proof}. This follows immediately from Theorem \ref{thm1} and Theorem \ref{thm3}.  \hfill $\square$

\medskip

To end this part, let us mention that   the entropy differential inequalities and the entropy power concavity inequalities in Theorem \ref{thm1}, Theorem \ref{thm2} and Theorem \ref{thm3} have been already proved in our 2017 preprint \cite{LL17}. In  \cite{LL17}, we have also extended the Entropy Power Concavity Inequality (EPCI)  for the Renyi entropy power to the porous medium equation associated with the Witten Laplacian on complete Riemannian manifolds with $CD(K, m)$-condition and on $(K, m)$-super Ricci flows.  
Due to the limit of space, we will put the nonlinear part into a forthcoming paper.  We would like also to mention that Wang and Zhang  \cite{Wang} proved the entropy concavity inequality for the $p$-Laplacian equation on compact Riemannian manifolds with non-negative Ricci curvature. Our work on EDI and EPCI  is independent of \cite{Wang}.

\medskip

\noindent{\bf Acknowledgement}. The second author would like to thank Prof. N. Mok for his suggestion which leads us to study Shannon entropy power  on manifolds and Ricci flow. We would like also to thank to
Prof.  Guangyue Han and Prof.  Fengyu Wang  for helpful discussions.

\medskip

\begin{flushleft}
\medskip\noindent

Songzi Li, School of Mathematics,  Renmin University of China,  Beijing, 100872, China

\medskip

Xiang-Dong Li, Academy of Mathematics and Systems Science, Chinese
Academy of Sciences, No. 55, Zhongguancun East Road, Beijing, 100190,  China,
E-mail: xdli@amt.ac.cn

and

School of Mathematical Sciences, University of Chinese Academy of Sciences, Beijing, 100049, China
\end{flushleft}

\end{document}